\documentclass[11pt,psamsfonts,twoside]{article}
\usepackage{amsmath,amssymb,latexsym}

\usepackage{diagrams}

\newcommand\qed{{\hspace*{\fill}Q.E.D.\vskip12pt plus 1pt}}

\newcommand\Pic[1]{\hbox{\rm Pic(}#1\hbox{\rm )}}

\newcommand\sI{{\cal I}}

\newcommand\sO{{\cal O}}

\newcommand\sX{{\cal M}}

\def\Coker{\operatorname{Coker}}
\def\Pic{\operatorname{Pic}}

\def\Ker{\operatorname{Ker}}
\def\id{\operatorname{id}}
\def\Im{\operatorname{Im}}

\newcommand\gra{\alpha}
\newcommand\grb{\beta}

\newcommand\comp{{\mathbb C}}

\newcommand\pn[1]{{\mathbb P}^{#1}}

\newcommand\proof{\noindent{\em Proof.}\ \ }

\newtheorem{theorem}{Theorem}[section]
\newtheorem{lemma}[theorem]{Lemma}
\newtheorem{corollary}[theorem]{Corollary}

\newtheorem{prop}[theorem]{Proposition}

\newtheorem{question}[theorem]{Question}

\newtheorem{definition}[theorem]{Definition}
\newtheorem{re}[theorem]{Remark}
\newtheorem{pargrph}[theorem]{}
\newtheorem{examp}[theorem]{Example}
\newtheorem{MM}[theorem]{ }
\newtheorem{res}[theorem]{Remarks}

\makeatletter
\ifnum\@ptsize=0  \fi
\ifnum\@ptsize=2  \fi
\renewcommand{\qed}{\hfill $\square$}

\newenvironment{rem*}{\begin{re}\em}{\end{re}}
\newenvironment{rems*}{\begin{res}\em}{\end{res}}
\newenvironment{example*}{\begin{examp}\em}{\end{examp}}
\newenvironment{definition*}{\begin{definition}\em}{\end{definition}}
\newenvironment{question*}{\begin{question}\em}{\end{question}}
\newenvironment{MM*}{\begin{MM}\em}{\end{MM}}

\newenvironment{prgrph*}[1]{\indent\begin{pargrph}{\bf #1.}\em\
}{\end{pargrph}}

\begin{document}

\title{Projective manifolds containing special curves}
\author{Lucian B\u{a}descu and Mauro C. Beltrametti\footnote{Universit\`a degli Studi di Genova, Dipartimento di Matematica, Via Dodecaneso 35, I-16146 Genova, Italy. e-mail:
badescu@dima.unige.it, beltrame@dima.unige.it}}

\date{ }

\maketitle

\begin{abstract} Let $Y$ be a smooth curve embedded in a complex projective manifold $X$ of dimension $n\geq 2$ with ample normal bundle $N_{Y|X}$. For every $p\geq 0$ let $\alpha_p$ denote the natural restriction maps $\Pic(X)\to\Pic(Y(p))$, where $Y(p)$ is the $p$-th infinitesimal neighbourhood of $Y$ in $X$. First one proves that for every $p\geq 1$ there is an isomorphism of abelian groups $\Coker(\gra_p)\cong\Coker(\gra_0)\oplus K_p(Y,X)$, where $K_p(Y,X)$ is a quotient of the $\mathbb C$-vector space $L_p(Y,X):=\bigoplus\limits_{i=1}^p H^1(Y, {\bf S}^i(N_{Y|X})^*)$ by a free subgroup of $L_p(Y,X)$ of rank strictly less than the Picard number of $X$. Then one shows that $L_1(Y,X)=0$ if and only if $Y\cong\mathbb P^1$ and $N_{Y|X}\cong\mathcal O_{\mathbb P^1}(1)^{\oplus n-1}$ (i.e. $Y$ is a quasi-line in the terminology of \cite{BBI}). The special curves in question are by definition those for which $\dim_{\mathbb C}L_1(Y,X)=1$. This equality is closely related with a beautiful classical result of B. Segre \cite{S}. It turns out that $Y$ is special if and only if either $Y\cong\mathbb P^1$ and $N_{Y|X}\cong\sO_{\pn 1}(2)\oplus\sO_{\pn 1}(1)^{\oplus n-2}$, or $Y$ is elliptic and $\deg(N_{Y|X})=1$. 
After proving some general results on manifolds of dimension $n\geq 2$ carrying special rational curves (e.g. they form a subclass of the class of rationally connected manifolds which is stable under small projective deformations), a complete birational classification of pairs $(X,Y)$ with $X$ surface and $Y$ special is given. Finally, one gives several examples of special rational curves in dimension $n\geq 3$.\end{abstract}

\section*{Introduction}

Let $X=\mathbb P^1\times\mathbb P^1$ be  a smooth quadric in $\pn 3$, and let $Y$ be a smooth curve of bidegree $(1,1)$ on $X$. Let $\Gamma$ be
a curve in $X$ of bidegree $(m,n)$ meeting transversely the conic $Y$ in $m+n$ distinct points $P_1,\ldots,P_{m+n}$. Let
$\gra_i$, $\grb_i$ be the two ruling lines of $X$ passing through $P_i$, let $\gamma_i$ be the tangent line of $\Gamma$ at $P_i$, and let $\theta_i$ be the tangent line of $Y$ at  $P_i$. These are four lines through $P_i$, contained in the projective tangent plane of $X$ at $P_i$. Thus it makes sense to consider the cross-ratios
$(\gra_i,\gamma_i,\theta_i,\grb_i)\in \comp$ of the four lines through the point $P_i$, $i=1,\ldots,m+n$. A result of B. Segre \cite{S}, \S 37 asserts that
\begin{equation}\label{Segre} \sum\limits_{i=1}^{m+n}(\gra_i,\gamma_i,\theta_i,\grb_i)=n.\end{equation} 
Conversely, given $m+n$ distinct points 
$P_1,\ldots,P_{m+n}\in Y$,  and a line $\gamma_i$ through each point $P_i$ contained in the tangent space of $X$ at $P_i$ satisfying \eqref{Segre}, then there exists a curve $\Gamma$ in $X$ of bidegree $(m,n)$ meeting $Y$ transversely only at the points $P_1,\ldots,P_{m+n}$ and such that $\gamma_i$ is the tangent line of $\Gamma$ at $P_i$, $i=1,\ldots,m+n$. In modern terminology this fact can be rephrased in terms of the Picard group of the first infinitesimal neighbourhood of $Y$ in $X$.

On the other hand, the classical condition of Reiss concerning the existence of a degree $d$ curve in $\pn 2$ intersecting a given line $Y$ in $d$ different prescribed points, with prescribed tangents and second-order conditions, can be reinterpreted in modern language in terms of the Picard group of the second infinitesimal neighbourhood of $Y$ in $\pn 2$ (see again \cite{S} and \cite{GH}, p. 698--699). 

These facts provide good motivation to study infinitesimal neighbourhoods of special curves in projective manifolds. More precisely, let $X$ be a complex projective manifold of dimension $n\geq 2$, and let $Y$ be a smooth connected curve of genus $g$ embedded in $X$ such that the normal bundle $N_{Y|X}$ of $Y$ in $X$ is ample. For every $p\geq 0$ we will denote by $Y(p)$ the $p$-th infinitesimal neighbourhood of $Y$ in $X$, i.e., $Y(p)$ is the algebraic scheme over $\comp$ whose underlying topological space coincides with the underlying topological space of $Y$, and whose structural sheaf $\sO_{Y(p)}$ is by definition $\sO_X/{\sI}^{p+1}$, where $\sI\subset \sO_X$ denotes the ideal sheaf of $Y$ in $X$. Of course $Y=Y(0)$. For every integer $p\geq 0$ we may consider the natural restriction maps
\begin{equation}\label{alpha}
\gra_p:{\rm Pic}(X)\to {\rm Pic}(Y(p)).\end{equation}
Then by Theorem \ref{Lemma1} below, for every $p\geq 1$ there exists an isomorphism
$$\Coker(\gra_p)\cong\Coker(\gra_0)\oplus K_p(Y,X),$$ 
where $K_p(Y,X)$ is a quotient of the $\mathbb C$-vector space 
$L_p(Y,X):=\bigoplus\limits_{i=1}^p H^1(Y, {\bf S}^i(N_{Y|X})^*)$ 
by a free subgroup of $L_p(Y,X)$ of rank $\leq\rho(X)-1$, where $\rho(X)$ is the Picard number of $X$.
Here ${\bf S}^i(E)$ denotes the $i$-th symmetric power of a vector bundle $E$.

The aim of this paper is to study the maps $\gra_p$, especially $\gra_1$ and $\gra_2$, when $L_p(Y,X)$ is of small dimension. For example, Theorem \ref{Theorem1} below describes the situation when $\dim_{\mathbb C}L_2(Y,X)$ is minimal. If $\dim(X)\leq 3$ or if $Y$ is not an elliptic curve this happens if and only if $Y\cong\mathbb P^1$ and $N_{Y|X}\cong\mathcal O_{\mathbb P^1}(1)^{\oplus n-1}$, i.e. if and only if
$Y$ is a quasi-line in $X$ in the terminology of \cite{BBI} (if $X$ is a surface this means that the embedding $Y\hookrightarrow X$ is Zariski equivalent to the embedding of a line in $\mathbb P^2$; this is the modern interpretation of Reiss' relation, see \cite{GH}, p. 698--699).

On the other hand, Corollary \ref{TheoremA} below asserts that
$L_1(Y,X)=0$ if and only if $\Coker(\alpha_1)$ is finite, or if and only if $Y$ is a quasi-line. 
Moreover, Theorem \ref{Theorem2} below takes care of the case 
$\dim_{\mathbb C}L_1(Y,X)=1$, in which case there are two possibilities: either
\begin{equation}\label{e:1}
Y\cong\mathbb P^1\;\;\text{and}\;\;N_{Y|X}\cong\sO_{\pn 1}(2)\oplus\sO_{\pn 1}(1)^{\oplus n-2},
\end{equation}
or $Y$ is an elliptic curve, $\deg(N_{Y|X})=1$ and the irregularity of $X$ is $\leq 1$.
Moreover, if we assume that $X$ is irregular and that $Y$ is $G3$ in $X$ (see Definition \ref{G3} below), then the canonical morphism of Albanese varieties 
${\rm Alb}(Y)=Y\to{\rm Alb}(X)$ is an isomorphism; in particular, the Albanese morphism $X\to{\rm Alb}(X)$ yields a retraction $\pi:X\to Y$ of the inclusion $Y\hookrightarrow X$.

In the case of surfaces one can say a lot more than Theorem \ref{Theorem2}. In fact, Theorem \ref{Theorem4} provides a very precise birational classification of pairs $(X,Y)$, with $X$ a smooth projective surface and $Y$ a smooth curve such that $(Y^2)>0$ and $\dim_{\mathbb C}L_1(Y,X)=1$.

The curves $Y$ in $X$ satisfying \eqref{e:1} are interesting from the point of view of varieties carrying quasi-lines (see
\cite{BBI} and \cite{IN}). Indeed by a result proved in \cite{IN} (Lemma 2.2), if $Y$ is such a curve and if $Z$ is a smooth
two-codimensional closed subvariety of $X$ meeting $Y$ at just one point transversely,
in the variety $X'$ obtained by blowing up $X$ along $Z$ the proper transform $Y'$ of $Y$ (via the blowing up morphism
$X'\to X$) becomes a quasi-line. In other words, any example of curves $Y$ in $X$ satisfying \eqref{e:1} provides
examples of projective manifolds containing quasi-lines. In section 4 we give several examples of projective manifolds $X$ carrying curves $Y$ satisfying \eqref{e:1}. One shows that the projective manifolds carrying curves satisfying \eqref{e:1} are rationally connected in the sense of \cite{KMM}, \cite{K}, and that they are stable under small projective deformations (Theorem \ref{Theorem3}).

Throughout this paper we shall use the standard terminology and notation in algebraic geometry.
All varieties considered are defined over the field $\mathbb C$ of complex numbers. 

\medskip

{\bf Acknowledgement.} The authors are grateful to the referee for some pertinent suggestions to improve the presentation and especially for pointing out an error in the proof of a previous formulation of Theorem \ref{Lemma1}.

\section{General results}\label{Background}\addtocounter{subsection}{1}\setcounter{theorem}{0}

Let $X$ be a complex projective manifold of dimension $n\geq 2$, and let $Y$ be a smooth connected curve of genus $g$ embedded in $X$ such that the normal bundle $N_{Y|X}$ of $Y$ in $X$ is ample.
For a non-negative integer $p$, we shall denote by $Y(p)$ the $p$-th infinitesimal 
neighbourhood $(Y,\sO_X/{\sI}^{p+1})$ of $Y$ in $X$ as in Introduction. Clearly $Y(0)=Y$.
Then for every $p\geq 1$ the truncated exponential sequence
$$0\to\sI^p/\sI^{p+1}\cong{\bf S}^p(N_{Y|X}^*)\to \sO_{Y(p)}^*\to\sO_{Y(p-1)}^*\to 0,$$
(in which $\sO_Z^*$ denotes the sheaf of multiplicative groups of nowhere vanishing functions on a scheme $Z$ and the first nontrivial map is the truncated exponential $u\mapsto 1+u$) yields the cohomology sequence
$$0\to H^0(Y,{\bf S}^p(N_{Y|X}^*))\to  H^0(Y(p),\sO_{Y(p)}^*)\to H^0(Y(p-1),\sO_{Y(p-1)}^*)\to$$
$$\to H^1(Y,{\bf S}^p(N_{Y|X}^*))\to {\rm Pic}(Y(p))\to{\rm Pic}(Y(p-1))\to H^2(Y,{\bf S}^p(N_{Y|X}^*))=0.$$
Since we work over a field of characteristic zero, ${\bf S}^p(N_{Y|X}^*)\cong{\bf S}^p(N_{Y|X})^*$ (see
\cite{Ha2}, Exercise 4.9, p. 114). Moreover, the hypothesis that $N_{Y|X}$ is ample implies that
${\bf S}^p(N_{Y|X})$ is also ample for every $p\geq 1$ (see \cite{Ha2}). From this it follows that 
$$H^0(Y,{\bf S}^p(N_{Y|X}^*))=0,$$
whence the map
$$a_p\colon H^0(Y(p),\sO_{Y(p)}^*)\to H^0(Y(p-1),\sO_{Y(p-1)}^*)$$ 
is an injective map of $\mathbb C$-algebras for every $p\geq 1$. But $H^0(Y(0),\sO_{Y(0)}^*)=H^0(Y,\sO_Y^*)=
\mathbb C^*=\mathbb C\setminus\{0\}$, and therefore $a_p$ is an isomorphism for every $p\geq 1$. It follows that the above cohomology sequence yields the exact sequence of abelian groups
\begin{equation}\label{1:p}
0\to H^1(Y,{\bf S}^p(N_{Y|X}^*))\to{\rm Pic}(Y(p))\to{\rm Pic}(Y(p-1))\to 0,\;\;\forall p\geq 1.\end{equation}
Since $H^1(Y,{\bf S}^p(N_{Y|X}^*))$ is a $\mathbb C$-vector space, the additive group $H^1(Y,{\bf S}^p(N_{Y|X}^*))$ is divisible 
(and hence injective), whence the exact sequence \eqref{1:p} splits for every $p\geq 1$. Then by induction we get
\begin{equation}\label{2:p}
{\rm Pic}(Y(p))\cong{\rm Pic}(Y)\oplus L_p(Y,X),\;\;\forall p\geq 1,
\end{equation}
where we put 
\begin{equation}\label{2':p}
L_p(Y,X):=\bigoplus\limits_{i=1}^pH^1(Y,{\bf S}^i(N_{Y|X}^*)).\end{equation}
Clearly, $L_p(Y,X)$ is a finite dimensional $\comp$-vector space.

\begin{theorem}\label{Lemma1} Let $X$ be a complex projective manifold of dimension $n\geq 2$, and let $Y$ be a smooth connected curve embedded in $X$ such that the normal bundle $N_{Y|X}$ of $Y$ in $X$ is ample. Then, for every $p\geq 1$, there exists an isomorphism
\begin{equation}\label{cok}\Coker(\gra_p)\cong\Coker(\gra_0)\oplus K_p(Y,X),\end{equation}
where $\gra_p\colon\Pic(X)\to\Pic(Y(p))$ is the map \eqref{alpha} and the abelian group $K_p(Y,X)$ is a quotient of the $\mathbb C$-vector space $L_p(Y,X)$ 
by a free subgroup of $L_p(Y,X)$ of rank $\leq \rho(X)-1$, with $\rho(X)$ the rank of the N\'eron-Severi group
of $X$ $($the Picard number of $X)$.\end{theorem}

\proof  Denote by $\beta_p\colon\Pic(Y(p))\to\Pic(Y)$ the natural restriction map and by $j\colon L_p(Y,X)\hookrightarrow\Pic(Y(p))$ the canonical inclusion into the direct sum (via the isomorphism \eqref{2:p}).
Now consider the commutative diagram
\begin{equation}\label{diagram}\begin{diagram}
&&0&&0&&0\\
&&\dTo&&\dTo&&\dTo\\
&&\Ker(\alpha_0)&\rTo^{\alpha'_p}&L_p(Y,X)&\rTo &K_p(Y,X)&\rTo &0\\
&&\dTo^i&&\dTo_j&&\dTo\\
&&\Pic(X)&\rTo^{\alpha_p}&\Pic(Y(p))&\rTo&\Coker(\alpha_p)&\rTo&0\\
&&\dTo&&\dTo_{\beta_p}&&\dTo_{\overline{\beta}_p}\\
0&\rTo&\Pic(X)/\Ker(\alpha_0)&\rTo^{\overline{\alpha}_0}&\Pic(Y)&\rTo&\Coker(\alpha_0)&\rTo&0\\
&&\dTo&&\dTo&&\dTo\\
&&0&&0&&0\\
\end{diagram}\end{equation}
in which the map $\overline{\alpha}_0$ is deduced from $\alpha_0$ by factorization, the map $\alpha'_p$ is induced by $\alpha_p$, the map $\overline{\beta}_p$ is induced by $\beta_p$ and $K_p(Y,X):=\Ker(\overline{\beta}_p)$.
Clearly, all columns and the second and the third rows are exact. Then by the snake lemma the first row is also exact. 
In particular, $K_p(Y,X)\cong L_p(Y,X)/\Im(\alpha'_p)$. Since $\mathbb Z$ is a principal ring and $L_p(Y,X)$ is an injective $\mathbb Z$-module, we infer that $K_p(Y,X)$ is also an injective $\mathbb Z$-module. This implies that the last column splits, which yields the isomorphism \eqref{cok}. Observe also that the subgroup $\Im(\alpha'_p)$ is torsion-free
since $L_p(Y,X)$ is a $\mathbb C$-vector space. Therefore $\Im(\alpha'_p)$ is free as soon as we know that $\Im(\alpha'_p)$ is a finitely generated group. Thus it remains to show that $\Im(\alpha'_p)$ is a finitely generated abelian group of rank $\leq\rho(X)-1$.

In view of decomposition \eqref{2:p}, the middle column splits, i.e. there exists a map $\eta\colon\Pic(Y(p))\to L_p(Y,X)$ such that $\eta\circ j=\id$. Clearly $\eta\circ\alpha_p\circ i=\eta\circ j\circ\alpha'_p=\alpha'_p$. Thus we get the map
$$\gamma_p:=\eta\circ\alpha_p\colon{\rm Pic}( X)\to L_p(Y,X),$$
such that $\gamma_p\circ i=\alpha'_p$. 

Observe now that the Picard scheme $\underline{{\rm Pic}}^0(X)$ is
an abelian variety  since  $X$ is smooth and projective (see \cite{FGA}, \'expos\'es 232, 236). Therefore 
$\gamma_p({\rm Pic}^0(X))=0$, since $L_p(Y,X)$ is an (additive) linear algebraic
group. Thus the map $\gamma_p$ factors through a map
$$\gamma'_p:{\rm NS}(X):={\rm Pic}(X)/{\rm Pic}^0(X)\to L_p(Y,X).$$
By the theorem of N\'eron-Severi, ${\rm NS}(X)$ is a finitely
generated abelian group of rank $\rho(X)\geq 1$ (since $X$ is projective).  
Since $L_p(Y,X)$ is a $\mathbb C$-vector space it follows that $\Im(\gamma_p)$ is a free abelian group of finite rank. Thus $\Im(\gamma_p)=\Im(\gamma'_p)$, and therefore also $\Im(\alpha'_p)$, is a finitely generated subgroup of $L_p(Y,X)$.
In fact, one can say more. Since $Y$ is a smooth projective curve, $\rho(Y)=1$. Note that the induced map ${\rm NS}(X)\to{\rm NS}(Y)$ is surjective after tensoring with $\mathbb Q$. Therefore the image of  $\Ker(\alpha_0)$ in ${\rm NS}(X)$ is  a (finitely generated) subgroup of rank equal to $\rho(X)-1$. Thus $\Im(\alpha'_p)$ is a free abelian subgroup of $L_p(Y,X)$ of rank $\leq\rho(X)-1$, which completes the proof of the theorem.\qed

\begin{rem*}\label{L} From Theorem \ref{Lemma1} it follows that $K_p(Y,X)=0$ if and only if $L_p(Y,X)=0$. Moreover,
if $\rho(X)=1$ we have $K_p(Y,X)\cong L_p(Y,X)$.
\end{rem*}

\begin{definition*}(\cite{BBI}) Let $Y$ be a smooth connected curve in the  projective manifold $X$ of dimension 
$n\geq 2$. The curve $Y$ is said to be a {\em quasi-line} in $X$ if $Y\cong\mathbb P^1$ and $N_{Y|X}\cong\mathcal O_{\mathbb P^1}(1)^{\oplus n-1}$. \end{definition*}

We are going to apply Theorem \ref{Lemma1} repeatedly. 
For instance we can easily compute the dimension of $L_2(Y,X)$
in the case when $Y\cong\pn 1$ and $N_{Y|X}\cong\sO_{\pn 1}(1)^{\oplus n-1}$, i.e. $Y$ is a quasi-line in $X$. Since $H^1(Y,N^*_{Y|X})\cong H^1(\pn 1,\sO_{\pn 1}(-1)^{\oplus n-1})=0$, we have in this case 
$$\dim_{\mathbb C} L_2(Y,X)=\dim_{\mathbb C} H^1(Y,{\bf S}^2(\sO_{\pn 1}(-1)^{\oplus n-1}))={{n}\choose{2}}=\frac{n(n-1)}{2}.$$

In particular, if $Y$ is a line in $X=\pn 2$, then the maps $\alpha_0$ and $\alpha_1$ are surjective and
$\dim_{\mathbb C} L_2(Y,X)=1$.
This is closely related to the so-called Reiss relation (see B. Segre \cite{S}, or also \cite{GH}, p. 698--699).

\begin{theorem}\label{Theorem1} Let $Y$ be a smooth connected curve embedded in a projective manifold $X$ of dimension $n\geq 2$
with normal bundle $N_{Y|X}$ ample. Then the following hold:\begin{enumerate}
\item[\em i)] If the curve $Y$ has genus  $g\neq 1$, then   $\dim_{\mathbb C} L_2(Y,X)\geq\dfrac{n(n-1)}{2}$. Moreover, the equality holds if
and only if
 $Y\cong\mathbb P^1$ and $N_{Y|X}\cong\sO_{\pn 1}(1)^{\oplus n-1}$ $($i.e., if and only if  $Y$ is a quasi-line in $X)$;
\item[\em ii)] If $Y$ is an elliptic curve, then  $\dim_{\mathbb C} L_2(Y,X)=(n+1)\deg(N_{Y|X})$.\end{enumerate}\end{theorem}

\proof By Theorem \ref{Lemma1} (with $p=2$) we have to compute $\dim_{\mathbb C}H^1(Y,{\bf S}^i(N_{Y|X}^*))$
for $i=1,2$. This follows from duality, Riemann-Roch, the fact that $H^0(Y,E^*)=0$ for every
ample vector bundle $E$ on $Y$, and the following formulae:
$$\deg({\bf S}^2(N_{Y|X}))=n\deg(N_{Y|X}), \;\;\text{and}\;\;{\rm rank}({\bf S}^2(N_{Y|X}))=
\frac{n(n-1)}{2}.$$
By the ampleness of $N_{Y|X}$ (which implies the fact that ${\bf S}^2(N_{Y|X})$ is also ample), and the standard formula
${\bf S}^2(N_{Y|X}^*)\cong{\bf S}^2(N_{Y|X})^*$, by duality we get 
$$H^1(Y,\omega_Y\otimes N_{Y|X})\cong H^0(Y,N_{Y|X}^*)=0,$$ and  $$H^1(Y,\omega_Y\otimes{\bf S}^2(N_{Y|X}))\cong
H^0(Y,{\bf S}^2(N_{Y|X})^*)=0.$$ 
Thus by duality and Riemann-Roch we have
\begin{multline}\label{3} \dim_{\mathbb C}H^1(Y,N_{Y|X}^*)=\dim_{\mathbb C}H^0(Y,\omega_Y\otimes N_{Y|X})=
\chi (Y,\omega_Y\otimes N_{Y|X})=\\
=\deg(N_{Y|X})+(n-1)(g-1),\end{multline}
and
\begin{multline}\label{4}\dim_{\mathbb C}H^1(Y,{\bf S}^2(N_{Y|X}^*))=\dim_{\mathbb C}H^0(Y,\omega_Y\otimes{\bf S}^2(N_{Y|X}))=\\
=\chi(Y,\omega_Y\otimes{\bf S}^2(N_{Y|X}))=n\deg(N_{Y|X})+\frac{n(n-1)}{2}(g-1).\end{multline}

From \eqref{3} and \eqref{4}  we get 
$$ \dim_{\mathbb C}L_2(Y,X)=(n+1)\deg(N_{Y|X})+\frac{(n+2)(n-1)}{2}(g-1).$$
Thus, if $g=1$, we get directly  ii).

Notice that since $N_{Y|X}$ is ample, $\deg(N_{Y|X})>0$. Thus if $g>1$ the last estimate yields  $\dim_{\mathbb C} L_2(Y,X)
\displaystyle>\frac{n(n-1)}{2}$. If instead $g=0$,  by a result of Grothendieck, we get
$$N_{Y|X}\cong\sO_{\pn 1}(a_1)\oplus\sO_{\pn 1}(a_2)\oplus\cdots\oplus
\sO_{\pn 1}(a_{n-1}),$$ with $ a_1\geq a_2\geq\cdots \geq a_{n-1}>0$
because $N_{Y|X}$ is ample. Then it is easily seen that the inequality of i) holds, with equality if and only if $a_1=\cdots=a_{n-1}=1$. This completes the proof of the theorem.\qed

\medskip

\begin{corollary} If in  Theorem {\em \ref{Theorem1}} we assume $n=2$ or $n=3$, then $\dim_{\mathbb C} L_2(Y,X)
\geq\dfrac{n(n-1)}{2}$, with equality if and only if
$Y$ is a quasi-line in $X$.\end{corollary}

\section{The first infinitesimal neighbourhood}\label{first}\addtocounter{subsection}{1}\setcounter{theorem}{0}

Now using Theorem \ref{Lemma1}, we proceed to analyze the map $\alpha_1$. As a direct consequence of Theorem \ref{Lemma1} we get the following result (see \cite{B1}, Theorem 14.2, which slightly improves Theorem (2.1) of \cite{BBI}; the latter  generalizes a result of d'Almeida \cite{dA} proved when $X$ is a surface and using different methods).

\begin{corollary}\label{TheoremA}  Let $Y$ be a smooth connected curve embedded in a projective
manifold $X$ of dimension $n\geq 2$
with normal bundle $N_{Y|X}$ ample. The following conditions are equivalent:
\begin{enumerate}
\item[\em i)] $L_1(Y,X)=0$.
\item[\em ii)] $K_1(Y,X)=0$.
\item[\em iii)] $\Coker(\alpha_1)$ is a finite group.
\item[\em iv)] $Y$ is a quasi-line.\end{enumerate}
Moreover, the map $\alpha_1$ is surjective if and only if $Y$ is a quasi-line and the map $\alpha_0$ is 
surjective.\end{corollary}

\proof  The equivalence i) $\Longleftrightarrow$ ii) follows from Remark \ref{L}. On the other hand, by duality we have $L_1(Y,X)=0$ if and only if $H^0(\omega_Y\otimes N_{Y|X})=0$. Using Riemann-Roch and the fact that every vector bundle on $\mathbb P^1$ is the direct sum of line bundles of the form $\mathcal O_{\mathbb P^1}(a)$, with $a\in\mathbb Z$, it is easy to see that the latter condition holds if and only if $Y\cong\pn 1$ and $N_{Y|X}\cong\sO_{\pn 1}(1)^{\oplus n-1}$. In particular, $L_1(Y,X)=0$ implies that $\Coker(\alpha_1)$ is finite because $\Pic(Y)\cong\mathbb Z$. Conversely, if $\Coker(\alpha_1)$ is finite then $L_1(Y,X)=0$ by Remark \ref{L}. The last statement is a direct consequence of decomposition \eqref{cok}. \qed

\medskip

To prove Theorem \ref{Theorem2} below we first need a definition and a lemma.

\begin{definition*}\label{G3} Let $Y$ be a closed subvariety of a projective irreducible variety $X$. We say that $Y$ is $G3$ in $X$
if the canonical map
$K(X)\to K(X_{/Y})$ is an isomorphism of rings, where $K(X)$ is the field of
rational functions of $X$, and $K(X_{/Y})$ is the ring of formal-rational functions of $X$ along $Y$ (see e.g. \cite{HM}, or also
\cite{Ha2}). In particular, if $Y$ is $G3$ in $X$, $K(X_{/Y})$ is a field.
We also say that $Y$ is $G2$ in $X$ if $K(X_{/Y})$ is a field and if
the above map makes $K(X_{/Y})$ a finite field extension of
$K(X)$.\end{definition*}

By a result of Hartshorne (see \cite{Ha2}, p. 198), if $X$ is smooth, $Y$ is connected and local complete intersection in $X$ and the normal bundle $N_{Y|X}$ is ample, then $Y$ is $G2$ in $X$. Thus, in our hypotheses from the beginning (i.e.
$Y$ is a smooth connected curve in the projective manifold $X$), $Y$ is always $G2$ in $X$. Moreover in the case when $X$ is a surface, $Y$ is $G3$ in $X$ whenever $N_{Y|X}$ is ample, i.e., $(Y^2)>0$. However, if $\dim X\geq 3$ and $Y\subset X$ is a curve with ample normal bundle, $Y$ is not necessarily  $G3$ in $X$ (see e.g. \cite{Ha2}, Exercise 4.10, p. 209, or also \cite{BBI}, Example (2.7)).

\begin{lemma}\label{Lemma2} Let $Y$ be an elliptic curve embedded in an irregular  projective manifold $X$ of dimension $n\geq 2$ with normal bundle $N_{Y|X}$ ample. Assume that $Y$ is $G3$ in $X$ $($this is always the case if $X$ is a surface$)$. Then the canonical morphism of Albanese varieties ${\rm Alb}(Y)=Y\to{\rm Alb}(X)$ $($induced by the inclusion $Y\hookrightarrow X)$ is an isomorphism. In particular, the Albanese morphism $f:X\to{\rm Alb}(X)$ yields a retraction $\pi:X\to Y$ of $Y\hookrightarrow X$.\end{lemma}

\proof By a general elementary result of Matsumura (see \cite{Ha2}, Exercise 4.15, p. 116), the morphism
$u:{\rm Alb}(Y)=Y\to{\rm Alb}(X)=:A$, associated to $Y\hookrightarrow X$, is surjective. In particular, $u$ is a finite \'etale morphism. Let $d$ be the degree of $u$. We have to prove that $d=1$.

Since $X$ is irregular and $Y$ is an elliptic curve we infer that $A$ is an elliptic curve, and since $f(X)$ generates $A$, the morphism $f:X\to A$ is surjective.  Consider now the cartesian diagram
$$\begin{diagram}
X':=X\times_A Y & \rTo^{u'}& X\\
\dTo^{f'} &   &  \dTo_f\\
Y &\rTo^u& A.\\
\end{diagram}$$
The inclusion $i:Y\hookrightarrow X$ yields a morphism $i':Y\to X'$ such that
$u'\circ i'=i$ and $f'\circ i'={\rm id}_Y$; in particular, $i'$ is a closed embedding and  $u'$ yields an isomorphism $i'(Y)\cong Y$. By a general elementary fact, $\dim(f(X))=\dim(A)=1$ implies that the morphism $f$ has connected fibers (see e.g. \cite{BS}, Lemma (2.4.5)). Since the above diagram is cartesian, it follows that $f'$  has also  connected fibers. Therefore $X'$ is connected (since $Y$ is so).

On the other hand, the morphism $u':X'\to X$ is finite and \'etale of degree $d$, because $u$ is so. Moreover, since $X$ is projective and nonsingular, $X'$ is also  projective and nonsingular. In other words, $X'$ is a projective manifold such that $u'$ and $i'$ define an \'etale neighbourhood of  $i:Y\hookrightarrow X$. In particular,  $u'$ yields an isomorphism of formal completions $\widehat{u'}:X'_{/i'(Y)}\cong X_{/Y}$, whence an isomorphism of rings of formal-rational functions ${\widehat{u'}}^*:K(X_{/Y})\cong
K(X'_{/i'(Y)})$.
Now, look at the commutative diagram
$$\begin{diagram}
K(X)& \rTo^{{u'}^*}& K(X')\\
\dTo &   & \dTo \\
K(X_{/Y}) &\rTo^{\widehat{u'}^*} &K(X'_{/i'(Y)}).\\
\end{diagram}$$
Since $Y$ is $G3$ in $X$ the first vertical map is an isomorphism. This and the isomorphism ${\widehat{u'}}^*$ imply 
that $\deg({u'}^*)=1$ (because the second vertical map is injective). But $\deg({u'}^*)=\deg(u')=d$, whence $d=1$. This completes the proof of the lemma.\qed

\medskip

Now we can prove the following result.

\begin{theorem}\label{Theorem2} Let $Y$ be a smooth connected curve of genus $g$ embedded in a  projective manifold $X$ of
dimension
$n\geq 2$ with normal bundle $N_{Y|X}$ ample.  Assume that $\dim_{\mathbb C} L_1(Y,X)=1$. Then $g\leq
1$. \begin{enumerate} 
\item[\em i)] If $g=0$ then $N_{Y|X}\cong\sO_{\pn 1}(2)\oplus\sO_{\pn 1}(1)^{\oplus n-2}$, and $X$ is rationally connected $($in the sense of {\rm \cite{KMM}}, cf. also {\rm \cite{K})}. 
\item[\em ii)] If $g=1$ then $\deg(N_{Y|X})=1$ and the irregularity of $X$ is $\leq 1$.
Assume moreover that $X$ is irregular and $Y$ is $G3$ in $X$. Then the canonical morphism of Albanese varieties 
${\rm Alb}(Y)=Y\to{\rm Alb}(X)$ is an isomorphism, and in particular, the Albanese morphism $X\to{\rm Alb}(X)$ yields a retraction $\pi:X\to Y$ of the inclusion $Y\hookrightarrow X$.\end{enumerate}\end{theorem}

\proof  The hypothesis says that $\dim_{\mathbb C}H^1(Y,N_{Y|X}^*)=1$. As in the proof of Theorem
\ref{Theorem1}, $\dim_{\mathbb C}H^1(N_{Y|X}^*)=\deg(N_{Y|X})+(n-1)(g-1)$, whence 
\begin{equation}\label{eq5} 1=\deg(N_{Y|X})+(n-1)(g-1).\end{equation}
Since $N_{Y|X}$ is ample, $\deg(N_{Y|X})\geq 1$, so that \eqref{eq5} implies $g\leq 1$.
Moreover, if $g=1$, it follows that $\deg(N_{Y|X})=1$.

i) If  $g=0$ then $Y\cong\mathbb P^1$, and \eqref{eq5} yields $\deg(N_{Y|X})=n$; moreover we get
$$N_{Y|X}=\sO_{\pn 1}(a_1)\oplus\sO_{\pn 1}(a_2)\oplus\cdots\oplus
\sO_{\pn 1}(a_{n-1}),\;\;\text{with}\;\; a_1\geq a_2\geq\cdots\geq a_{n-1}.$$
It follows that $\deg(N_{Y|X})=a_1+a_2+\cdots+a_{n-1}$, and since $N_{Y|X}$ is ample, $a_{n-1}>0$.
Since $\deg(N_{Y|X})=n$ we get $a_1+a_2+\cdots+a_{n-1}=n$, whence $a_1=2$ and $a_2=\cdots=a_{n-1}=1$.
Then it is a general fact that $X$ is rationally connected (see  \cite{K}, \cite{KMM}). 

ii) When $g=1$, the ampleness of $N_{Y|X}$ and the result of Matsumura quoted in the proof of Lemma \ref{Lemma2} imply that the map ${\rm Alb}(Y)=Y\to{\rm Alb}(X)$ is surjective, whence the irregularity of $X$ is $\leq 1$. Finally, if the irregularity of $X$ is $1$, Lemma \ref{Lemma2}  proves the remaining part of ii), thus completing the proof of  the theorem.\qed

\medskip

\begin{rems*}\label{Remarks}  i) The hypothesis in Theorem \ref{Theorem2}, ii) that $Y$ is $G3$ in $X$ is not very restrictive. Indeed, as we remarked above, the ampleness of $N_{Y|X}$ implies by the result of Hartshorne quoted above that $Y$ is in any case $G2$ in $X$. Then by a result of Hartshorne-Gieseker (see \cite{Gi1}, Theorem 4.3) there is a finite morphism $f:X'\to X$ with the following properties: the inclusion $Y\hookrightarrow X$ lifts to a closed embedding $j:Y\hookrightarrow X'$ such that $f$ is \'etale along $j(Y)$ (i.e. $(X',j(Y))$ is an \'etale neighbourhood of $(X,Y)$) and $j(Y)$  is $G3$ in $X'$. In particular, $X'$ is nonsingular along $Y$ and $N_{j(Y)|X'}\cong N_{Y|X}$ is ample. Desingularizing $X'$ away $j(Y)$ we get even a projective manifold $\widetilde{ X}$ containing $Y$
such that $N_{Y|\widetilde{X}}\cong N_{Y|X}$ is ample and $Y$ is $G3$ in $\widetilde{X}$. Moreover, if $X$ is irregular, $\widetilde{X}$ is also irregular (both having irregularity $1$ since $N_{j(Y)|\widetilde{X}}\cong N_{Y|X}$ is ample and the morphism $\widetilde{X}\to X$ is surjective). 

\medskip

ii) Let $Y\subset X$ be as in Theorem \ref{Theorem2}, and assume that the irregularity of $X$ is $1$ (in particular, $Y$ is an elliptic curve). We claim that the normal exact sequence 
$$0\to T_Y=\sO_Y\to T_X|Y\to N_{Y|X}\to 0$$
splits. Indeed, if $Y$ is $G3$ in $X$, then the retraction $X\to Y$ of $Y\hookrightarrow X$ yields the desired splitting. Otherwise, use the previous remark to lift the embedding
$Y\hookrightarrow X$ to $j:Y\hookrightarrow X'$ such that $j(Y)$ is $G3$ in $X'$ and $f:X'\to X$  is \'etale along $j(Y)$. Then by Lemma \ref{Lemma2}
there exists a retraction $X'\to Y$ for $j$, so that the normal exact sequence  
$$0\to T_Y=\sO_Y\to T_{X'}|Y\to N_{Y|X'}=N_{Y|X}\to 0$$
splits. Since $f$ is \'etale along $j(Y)$ then the splitting of the latter normal sequence implies the splitting of the former normal sequence.

\medskip

iii) Assume, as in Theorem \ref{Theorem2}, that $Y$ is a smooth connected curve of genus $g$ embedded in a  projective manifold $X$ of dimension $n\geq 2$ with normal bundle $N_{Y|X}$ ample. Then the arguments of the proof of Theorem \ref{Theorem2}, i) yield in fact the following more general statement: if  $L_1(Y,X)$ is of dimension $h<n$
(respectively $h=n$) then $g\leq 1$ (respectively $g\leq 2$, and if $g=2$ then $\deg(N_{Y|X})=1$).
Indeed, instead of equality \eqref{eq5} we have
$$h=\deg(N_{Y|X})+(n-1)(g-1).$$
If $h<n$, since $\deg(N_{Y|X})\geq 1$ we cannot have $g\geq 2$. If $h=n$ then $g\leq 2$,  with $\deg(N_{Y|X})=1$ if
$g=2$. 

If $g=1$ then $\deg(N_{Y|X})=h$ and the irregularity of $X$ is $\leq 1$.
If instead $g=0$ one has $\deg(N_{Y|X})=h+n-1$ and
$$N_{Y|X}\cong\sO_{\pn 1}(a_1)\oplus\cdots\oplus\sO_{\pn 1}(a_h)\oplus
\sO_{\pn 1}(1)^{\oplus n-h-1},$$
with $a_1\geq a_2\geq\cdots\geq a_h\geq 1$ and $\sum\limits_{i=1}^ha_i=2h$. In particular, $X$ is rationally connected.
\end{rems*}

Now we recall the following result of Ionescu-Naie \cite{IN}, Lemma (2.2): 

\begin{theorem}[\cite{IN}]\label{quasi} Let $X$ be a projective manifold of dimension $n\geq 2$ and let $Y$ be a smooth rational curve in $X$ with normal bundle of the form
$$N_{Y|X}\cong\sO_{\pn 1}(a_1)\oplus\sO_{\pn 1}(a_2)\oplus\cdots\oplus\sO_{\pn 1}(a_{n-1})\;\;\text{with}\;\; a_1\geq a_2\geq\cdots\geq a_{n-1}.$$ 
Let $Z\subset X$ be a general smooth $2$-codimensional subvariety of $X$ meeting $Y$ transversely in one point. Let $f:\widetilde{X}\to X$ be the blowing up of $X$ along $Z$ and let $\widetilde{Y}$ be the proper transform of $Y$ via $f$. Then the normal bundle of $\widetilde{Y}$ in $\widetilde{X}$ is given by
$$N_{\widetilde{Y}|\widetilde{X}}\cong\sO_{\pn 1}(a_1-1)\oplus\sO_{\pn 1}(a_2)\oplus\cdots\oplus\sO_{\pn 1}(a_{n-1}).$$ \end{theorem}

\begin{corollary}\label{quas1} Let $X$ be a projective manifold of dimension $n\geq 2$ and let $Y$ be a smooth rational curve in $X$ with normal bundle of the form $N_{Y|X}\cong\sO_{\pn 1}(2)\oplus\sO_{\pn 1}(1)^{\oplus n-2}$. Let $Z\subset X$ be a general smooth $2$-codimensional subvariety of $X$ meeting $Y$ transversely in one point. Let $f:\widetilde{X}\to X$ be the blowing up of $X$ along $Z$ and let $\widetilde{Y}$ be the proper transform of $Y$ via $f$. Then $\widetilde{Y}$ is a quasi-line in $\widetilde{X}$.\end{corollary}

\begin{rem*}\label{quas2} Let $Y$ be a smooth rational curve in the projective $n$-fold $X$ (with $n\geq 2$) such that 
$N_{Y|X}\cong\sO_{\pn 1}(2)\oplus\sO_{\pn 1}(1)^{\oplus n-2}$. Since by Bertini there always exists smooth $2$-codimensional subvarieties $Z$ of $X$ meeting $Y$ transversely in one point, Corollary \ref{quas1} shows that as soon as we start with such a pair $(X,Y)$ we easily produce a projective $n$-fold $\widetilde{X}$ (dominating $X$) and a quasi-line $\widetilde{Y}$ in $\widetilde{X}$.\end{rem*}

As proved in \cite{KMM}, the rationally connected manifolds are stable under any projective deformation. The next result shows that the projective manifolds of dimension $n$ carrying nonsingular rational curves with normal bundle $\sO_{\pn 1}(2)\oplus\sO_{\pn 1}(1)^{\oplus n-2}$ are stable under small (but not under global) projective deformations. A similar result holds for quasi-lines, see \cite{BBI}, (3.10). 

\begin{theorem}\label{Theorem3} Any small projective deformation of a projective manifold $X$ 
of dimension $n\geq 2$ containing a smooth rational curve $Y$ such that $N_{Y|X}\cong
\sO_{\pn 1}(2)\oplus\sO_{\pn 1}(1)^{\oplus n-2}$ is a projective manifold containing a
smooth rational curve with normal bundle isomorphic to
$\sO_{\pn 1}(2)\oplus\sO_{\pn 1}(1)^{\oplus n-2}$.\end{theorem}

\proof Let $f\colon\sX\to T$ be a smooth projective morphism such that there is a point $t_0\in T$ with the property that $f^{-1}(t_0)\cong X$. By
taking an appropriate base change we may assume that $T$ is a smooth curve. We may view $Y$ as a curve in $\sX$. 
The proof of the openness of the deformations of rationally
connected manifolds works in our situation as well (see the first part of
the proof of Proposition (2.13) of \cite{MP}, p. 107). In fact in our case it becomes
even simpler, working with the Hilbert scheme (instead of the ${\rm Hom}$-scheme).
In fact, consider the canonical exact sequence 
\begin{equation} 0\to N_{Y|X}\to N_{Y|\sX}\to N_{X|\sX}|Y\cong\sO_Y\to 0.
\end{equation}
Since $H^1(Y,N_{Y|X})\cong H^1(\pn 1,\sO_{\pn 1}(2)\oplus\sO_{\pn 1}(1)^{\oplus n-2})=0$,
the exact  sequence  splits. Therefore we get 
$$H^0(Y,N_{Y|\sX})\cong H^0(Y,N_{Y|X})\oplus\mathbb C\;\;\text{and}\;\; H^1(Y,N_{Y|\sX})=0.$$
It follows that there exists an one parameter family of curves $\{Y_s\}_{s\in D}$ (parametrized by the unit disk $D$) such that $Y_0=Y$ and $Y_s$ is not contained in  $X_{t_0}=X$ for $s\neq 0$. Since $Y_0$ is contained in the fiber  $f^{-1}(t_0)$, for each $s\in D$ the curve $Y_s$ is contained in some fiber $X_{t_s}=f^{-1}(t_s)$, and the morphism defined by $s\mapsto t_s$ is unramified near $0$. Clearly $N_{Y|X}=\sO_{\pn 1}(2)\oplus\sO_{\pn 1}(1)^{\oplus n-2}$ is ample. Since ampleness is an open condition and $\deg(N_{Y_s|X_{t_s}})=\deg(N_{Y|X})=n$ for all $s\in T$, it
follows that $N_{Y_s|X_{t_s}}\cong\sO_{\pn 1}(2)\oplus\sO_{\pn 1}(1)^{\oplus n-2}$ for
$s$ near $0$. Finally, since $T$ is a smooth curve, there is an open neighbourhood of $t_0$ in $T$ over which $X_t$ contains a smooth rational curve $Y_t$ with $N_{Y_t|X_t}\cong
\sO_{\pn 1}(2)\oplus\sO_{\pn 1}(1)^{\oplus n-2}$.\qed

\begin{rem*}\label{Hilbert} Let $X$ be an
$n$-dimensional manifold $X$ containing a smooth rational curve $Y$  with ample
normal bundle such that $\dim_\comp L_1(Y,X)=1$. Then by Theorem
\ref{Theorem2}, i), we have $N_{Y|X}\cong \sO_{\pn 1}(2)\oplus\sO_{\pn 1}(1)^{\oplus (n-2)}$. Let ${\rm Hilb}_Y(X)$ be the Hilbert scheme of $Y$. Then standard considerations yield the following facts: ${\rm Hilb}_Y(X)$ is smooth at the point corresponding to $Y$, the general embedded deformations of $Y$ are smooth rational curves having the same normal bundle as $Y$, and their union is dense in $X$. Moreover, through any three general points of $Y$ there pass only finitely many smooth rational curves from the given
family.\end{rem*}

\section{The case of surfaces}\label{Surfaces}\addtocounter{subsection}{1}\setcounter{theorem}{0}

Now we want to look more closely to what happens in the case when $n=2$, i.e. when $X$ is a surface. First let us give some examples.

\begin{example*}\label{Example1} Let $X:=\pn 1\times\mathbb P^1$ and let $Y\in|\sO_
{\pn 1\times\mathbb P^1}(1,1)|$ be any smooth curve. Then $Y\cong\mathbb P^1$,
$N_{Y|X}\cong\sO_{\pn 1}(2)$, whence $\dim_{\mathbb C}H^1(Y,N_{Y|X}^*)=1$. By Theorem \ref{Lemma1}, $\Coker(\alpha_1)\cong\mathbb C/F$, with $F$ a free subgroup of $(\mathbb C,+)$ of rank $\leq 1$. This was the case classically studied by B. Segre (see \cite{S}, \S 37).\end{example*}

\begin{example*}\label{Example2} Let 
$X:=\mathbb P(\sO_{\pn 1}\oplus\sO_{\pn 1}(-2))$ be the Segre-Hirzebruch surface
$\mathbb F_2$. Let $C_0$ be the minimal section of the canonical projection $\pi\colon X\to\mathbb P^1$
($(C_0^2)=-2$), and let $Y$ be a section of $\pi$ such that $(Y^2)=2$ and $Y\cap C_0=
\varnothing$. Clearly, $Y\cong\mathbb P^1$, $N_{Y|X}\cong\sO_{\pn 1}(2)$, and since
$Y$ is a section of $\pi$ the map $\alpha_0$ is surjective. \end{example*}

Note that  Examples \ref{Example1} and \ref{Example2} of embeddings of $\pn 1$ into $\pn 1\times\mathbb P^1$
and $\mathbb F_2$ respectively are not Zariski equivalent. Indeed, if we blow down the minimal
section $C_0$ of $\mathbb F_2$ we get the projective cone $V$ in $\mathbb P^3$ over the conic of $\mathbb P^2$ of equation $x_1^2=x_0x_2$. Then the conclusion follows from the facts that the image of $Y$ in $V$ is an ample divisor on $V$ and, on the other hand, in Example \ref{Example1} the curve $Y$ is ample on $\pn 1\times\mathbb P^1$.

A result of Gieseker (see \cite{Gi1}, Theorem 4.5) together with the fact that the embeddings of 
$\pn 1$ of Examples \ref{Example1}, \ref{Example2} are not Zariski equivalent implies that they are not formally equivalent either.
More precisely, in both cases we  have $(Y^2)>0$, whence by \cite{HM} or by \cite{Ha2},
$Y$ is $G3$ in $X$. Then the above-mentioned result of Gieseker tells us that if the two
embeddings were formally equivalent, then they would be also Zariski equivalent.

\begin{example*}\label{Example3} Let $B$ be an elliptic curve and $L$ a line bundle of degree
one on $B$. Since $H^1(B,L^{-1})\neq 0$ (in fact, $\dim_{\mathbb C}H^1(B,L^{-1})=1$), 
there exists a non splitting exact sequence of vector bundles
$$0\to\sO_B\to E\to L\to 0.$$
Then a result of Gieseker \cite{Gi2} shows that $E$ is ample because $L$ is so.
Put $X:=\mathbb P(E)$, and let $\pi\colon X\to B$ be the canonical projection. Let also $Y$ be the 
section of $\pi$ that corresponds to the surjection $E\to L$. Then $Y\in|\sO_{\mathbb P(E)}(1)|$, and in particular, $Y$ is an ample Cartier divisor on $X$ with normal bundle $N_{Y|X}=L$. Clearly, the 
map $\alpha_0$ is surjective, $(Y^2)=\deg(E)=1$ and $N_{Y|X}\cong L$. It follows that $\dim_{\mathbb C}H^1(Y,N_{Y|X}^*)=1$. \end{example*}

\begin{example*}\label{Example4} Let $B$ be an elliptic curve and $L$ a line bundle of degree one on $B$. Set $F:=\sO_B\oplus L$ and $X:=\mathbb P(F)$. Let $Y$ be the section of the canonical projection $\pi\colon X\to B$ corresponding 
to the canonical map $F\to L$. Then again
$N_{Y|X}\cong L$. There is another section $C_0$ of $\pi$ (the minimal section corresponding to the map $F\to\sO_B$) such that $(C_0^2)=-1$
 and $C_0\cap Y=\varnothing$. Then $C_0$ can
be blown down to get the projective cone $X'$ over the polarized curve $(B,L)$, i.e.
$X'\cong{\rm Proj}(S[T])$, where $S:=\bigoplus\limits_{i=0}^{\infty}H^0(B,L^i)$, $T$ is an indeterminate
over $S$ and the grading of $S[T]$ is given by $\deg(sT^j)=i+j$ whenever $s\in S$ is a
homogeneous element of degree $i$. Let $f\colon X\to X'$ be the blowing down morphism and set $Y':=f(Y)$. Then $Y'$ is an elliptic curve (isomorphic to $B$) such that $Y'$ is embedded in the smooth locus of $X'$ with normal bundle isomorphic to $L$.\end{example*}

\medskip

Consider Examples \ref{Example3} and \ref{Example4} with the same $B$ and $L$. Then $Y\cong Y'$ and $N_{Y|X}\cong N_{Y'|X'}$. On the other hand, exactly as in Examples \ref{Example1} and \ref{Example2}, one shows that these two embeddings are not formally equivalent (and hence not Zariski equivalent either).

To draw a consequence of Theorem \ref{Theorem2} we need to recall two well known results.

\begin{theorem}[\cite{Giz}]\label{TheoremB}  Let $X$ be a normal projective surface containing $Y=\pn 1$ as an ample Cartier divisor. Then, up to isomorphism, one has one of the following cases:
\begin{enumerate}
\item[\em i)] $X=\mathbb P^2$ and $Y$ is either a line or a conic; or
\item[\em ii)] $X=\mathbb F_e=\mathbb P(\sO_{\pn 1}\oplus\sO_{\pn 1}(-e))$ and $Y$ is a section of the canonical projection $\pi\colon\mathbb F_e\to\mathbb P^1$; or
\item[\em iii)] $X$ is the projective cone in $\mathbb P^{s+1}$ over the rational normal curve of 
degree $s$ in $\mathbb P^s$, and $Y$ is the intersection of $X$ with the hyperplane at infinity.\end{enumerate}
\end{theorem}

Theorem \ref{TheoremB} is classical, a modern reference for it is \cite{Giz}.

\begin{theorem}[\cite{F},\cite{B2}]\label{TheoremC} Let $X$ be a normal projective  surface containing an elliptic curve $Y$ as an ample Cartier divisor. Then one has one of the following cases:
\begin{enumerate}
\item[\em i)] $X$ is a $($possibly singular$)$ Del Pezzo surface {\rm (}i.e. a rational surface with at most rational double points as
singularities and with ample anticanonical class{\rm )}, and $-Y$ is a canonical divisor of $X$; or
\item[\em ii)] There exists an elliptic curve $B$ and an ample rank two vector bundle $E$
on $B$ such that $X\cong{\mathbb P}(E)$ and $Y\in|\sO_{\mathbb P(E)}(1)|$ $($in particular, $Y$
is a section of the canonical projection ${\mathbb P}(E)\to B)$; or
\item[\em iii)] $X$ is the projective cone over the polarized curve $(Y,N_{Y|X})$ {\rm (}i.e.
$X\cong{\rm Proj}(S[T])$, where $S=\bigoplus\limits_{i=0}^{\infty} H^0(Y,N_{Y|X}^{\otimes i})$, $T$ is an indeterminate over $S$ and the 
gradation of $S[T]$
is given by $\deg(sT^j)=\deg(s)+j$, whenever $s$ is a homogeneous element of $S${\rm )} and $Y$ is embedded in $X$ as the infinite section.
\end{enumerate}
\end{theorem}

Theorem \ref{TheoremC} is a generalization of a classical result, see \cite{F} if $X$ is smooth, and \cite{B2}, p. 3, if $X$ is singular.

Now we can prove the main result of this section:

\begin{theorem}\label{Theorem4} Let $X$ be a smooth projective surface and $Y$ a smooth connected  curve on $X$ such that $(Y^2)>0$ and $\dim_{\mathbb C}L_1(Y,X)=1$. Then there exists a birational morphism $\varphi\colon X\to X'$ and a Zariski open neighbourhood $U$  of $Y$ in $X$ such that the restriction $\varphi|U\colon U\to\varphi(U)$ is a biregular isomorphism,  $Y':=\varphi(Y)$ is an ample Cartier divisor on $X'$, and $(X',Y')$ is one of the following pairs:
\begin{enumerate}
\item[\em i)] $X'\cong\mathbb F_0=\pn 1\times\mathbb P^1$ and $Y'\in|\sO(1,1)|$; or
\item[\em ii)] $X'$ is isomorphic to the quadratic normal cone in $\mathbb P^3$ of equation $x_1^2=x_0x_2$, and $Y'$ is the intersection of $X'$
with the hyperplane $x_3=0$; or
\item[\em iii)] $Y$ is an elliptic curve, $(Y^2)=1$, and there exists an exact sequence of vector bundles on $Y$
$$0\to\sO_Y\to E\to N_{Y|X}\to 0$$
with $E$ ample such that $X'\cong\mathbb P(E)$ and $Y'\in|\sO_{\mathbb P(E)}(1)|$; or
\item[\em iv)] $Y$ is an elliptic curve such that $(Y^2)=1$ and $X'$ is the projective cone over the polarized curve $(Y,N_{Y|X})$ $($i.e.
$X\cong{\rm Proj}(S[T])$, where $S=\bigoplus\limits_{i=0}^{\infty} H^0(Y,N_{Y|X}^{\otimes i})$, $T$ is an indeterminate over $S$ and the gradation of $S[T]$
is given by $\deg(sT^j)=\deg(s)+j$, whenever $s$ is a homogeneous element of $S)$ and $Y'$ is embedded in $X'$ as the infinite section
$($i.e. $Y'=D_+(T))$; or
\item[\em v)] $Y$ is an elliptic curve such that $(Y^2)_{X'}=1$ and $X'$ is a $($possibly singular$)$ Del Pezzo surface of degree $1$ and
$-Y$ is a canonical divisor of $X'$.
$($These surfaces are classified in {\rm \cite{De}.)}\end{enumerate}
\end{theorem}
\proof If $X$ is a surface then $Y$ is a divisor on $X$. Since $N_{Y|X}$ is ample, by \cite{Ha2}, Theorem  4.2, p. 110,  there exists a birational isomorphism $\varphi\colon X\to X'$ with the following properties:

$-$ $X'$ is a normal projective surface,

$-$ there is a Zariski open neighbourhood $U$ of $Y$ in $X$ such that the restriction
$\varphi|U\colon U\to\varphi(U)$ is a biregular isomorphism, and

$-$ $Y':=\varphi(Y)$ is an ample Cartier divisor on $X'$.

Note that in loc. cit. one first proves that the linear system $|mY|$ is base point free for $m\gg 0$. Then $\varphi$ is gotten from 
the morphism associated to $|mY|$, for $m\gg 0$, by passing to the Stein factorization.

Now, by Theorem \ref{Theorem2}, $g\leq 1$; moreover, $(Y^2)=2$ if $g=0$, and $(Y^2)=1$ if $g=1$. Now the 
classification of the normal projective surfaces $X'$ supporting a smooth rational or a smooth elliptic curve $Y'$ as an ample Cartier 
divisor is given by Theorems \ref{TheoremB} and \ref{TheoremC} above. 

If $g=0$, $({Y'}^2)_{X'}=(Y^2)_X=2$, and then we apply Theorem \ref{TheoremB}. In  case i) of \ref{TheoremB} we have
that $({Y'}^2)_{X'}$ is $1$ or $4$, whence this case is ruled out. Moreover, $({Y'}^2)_{X'}=2$ can be realized in cases ii) or iii) of Theorem \ref{TheoremB}  either if $X'\cong\mathbb F_0=\pn 1\times\mathbb P^1$ and $Y'\in|\sO(1,1)|$ (and this corresponds to Example \ref{Example1} above), or if $Y'$ is isomorphic to the quadratic normal cone in $\mathbb P^3$ of equation $x_1^2=x_0x_2$ (which corresponds to Example \ref{Example2} above). 

If $g=1$ and the surface $X$ is rational, by Theorem \ref{TheoremC}, i), $X'$ is a (possibly singular)
Del Pezzo surface and $-Y$ is a canonical divisor of $X'$. Moreover, since $(Y^2)_{X'}=1$,
the degree of $X'$ is $1$. 

Assume now $X$ not rational. Then $X'$ is a surface as in each of cases ii) or iii) of Theorem \ref{TheoremC}.
In both cases, by Theorem \ref{Theorem2}, we have $\deg(N_{Y|X})=1$, whence $({Y'}^2)_{X'}=1$.

If we are in case ii) (of Theorem \ref{TheoremC}), then $X'\cong\mathbb P(E)$, with $E$ an ample rank two vector bundle over
an elliptic curve $B$, and $Y'\in|\sO_{\mathbb P(E)}(1)|$. Let $\pi\colon\mathbb P(E)\to B$ be the canonical projection. Then the exact sequence
$$0\to\sO_{\mathbb P(E)}\to\sO_{\mathbb P(E)}(1)\cong\sO_{X}(Y)\to N_{Y'|X'}\cong N_{Y|X}\to 0$$
yields the cohomology sequence
$$0\to\pi_*(\sO_{\mathbb P(E)})\to\pi_*(\sO_{\mathbb P(E)}(1))\to\pi_*(N_{Y'|X'})\to
R^1\pi_*(\sO_{\mathbb P(E)})=0,$$
or else,
$$0\to \sO_Y\to E\to N_{Y|X}\to 0.$$ So we get case iii) of our statement.

Finally, case iii) of Theorem \ref{TheoremC} yields case iv).\qed

\begin{rem*} In Theorem \ref{Theorem4}, the cokernel of the restriction map $\alpha_0\colon{\rm Pic}(X)\to{\rm Pic}(Y)$ is finite (and in fact the map $\gra_0$
is surjective) if and only if $(X,Y)$ is in one of cases i)--iv) (see \cite{B1}, \S 14).\end{rem*}

\section{Examples in higher dimension}\label{Examples}\addtocounter{subsection}{1}\setcounter{theorem}{0}

To give the first examples of curves $Y\cong\mathbb P^1$ on projective $n$-folds $X$, such that $n\geq 3$ and
$N_{Y|X}\cong\sO_{\pn 1}(2)\oplus\sO_{\pn 1}(1)^{\oplus n-2}$, we need the following simple lemma.

\begin{lemma}\label{Lemma3} Let $X$ be a smooth projective variety of dimension $n\geq 3$, and let 
$X'$ be a smooth irreducible hypersurface of $X$ such that:
\begin{enumerate}
\item[{\rm i)}] $X'$ contains a quasi-line $Y$, i.e. there is a curve $Y$ on $X'$ such that $Y\cong\mathbb P^1$ and $N_{Y|X'}\cong\sO_{\pn 1}(1)^{\oplus n-2}$; and
\item[{\rm ii)}]  $(N_{X'|X}\cdot Y)=2$.\end{enumerate}
Then $N_{Y|X}\cong\sO_{\pn 1}(2)\oplus\sO_{\pn 1}(1)^{\oplus n-2}$.\end{lemma}

\proof Consider the canonical exact sequence of normal bundles
$$0\to N_{Y|X'}\to N_{Y|X}\to N_{X'|X}|Y\to 0, $$
in which $N_{Y|X'}\cong\sO_{\pn 1}(1)^{\oplus n-2}$  and $N_{X'|X}|Y\cong\sO_{\pn 1}(2)$ by conditions i) and ii) respectively. Since $H^1(\pn
1,\sO_{\pn 1}(-1)^{\oplus n-2})=0$, this exact sequence splits to give the conclusion.\qed

\begin{example*}\label{Example5} Let $v_2\colon\mathbb P^{n-1}\hookrightarrow\mathbb P^m$ be the $2$-fold Veronese embedding of $\mathbb P^{n-1}$ 
(with $m=(n^2+n-2)/2$). In Lemma \ref{Lemma3} we take $X'=v_2(\mathbb P^{n-1})$ and $Y=v_2(L)$, with $L$ a line in $\mathbb P^{n-1}$. Then $Y$ is
a  smooth conic in $\mathbb P^m$, contained in $X'$. In particular, $N_{Y|X'}\cong\sO_{\pn 1}(1)^{\oplus n-2}$. Let
 now $Z\subset\mathbb P^{m+1}$ be the cone over $X'$ with vertex a point $z\in\mathbb P^{m+1}\setminus\mathbb P^m$. Notice that $Z$ is isomorphic to the
weighted projective space
$\mathbb P^n(1,\ldots,1,2)$. Then $Z$ contains $X'\cong\mathbb P^{n-1}$ as an ample Cartier
divisor such that $N_{X'|Z}\cong\sO_{\mathbb P^{n-1}}(2)$. Let $X$ be the blowing up of $Z$ at $z$. Then
 $X\cong\mathbb P(\sO_{\mathbb P^{n-1}}(2)\oplus\sO_{\mathbb P^{n-1}})$, $X$ still 
contains $X'$, and $N_{X'|X}\cong N_{X'|Z}\cong\sO_{\mathbb P^{n-1}}(2)$, whence 
$(N_{X'|X}\cdot Y)=2$. Therefore by Lemma \ref{Lemma3} we get
\begin{equation}\label{eq7} 
N_{Y|X}\cong\sO_{\pn 1}(2)\oplus\sO_{\pn 1}(1)^{\oplus n-2}.\end{equation}

We call this example {\em the standard example} of a smooth rational curve $Y$ in an $n$-fold
$X$ satisfying \eqref{eq7}. Example \ref{Example5} is a higher dimensional analogue of Example \ref{Example2}. Note that $X$ is also isomorphic to the projective closure $\mathbb P(\sO_{\mathbb P^{n-1}}(-2)\oplus\sO_{\mathbb P^{n-1}})$ of the geometric vector bundle $\mathbb V(\sO_{\mathbb P^{n-1}}(-2))$.\end{example*}

\begin{example*}\label{Example6} Consider the projective bundle $X':=\mathbb P(T_{\mathbb P^d})$ associated to the tangent bundle $T_{\mathbb P^d}$ of 
$\mathbb P^d$, with $d\geq 2$. In particular, $\dim(X')=2d-1$. It is well known that $X'$ contains
quasi-lines $Y$ (see e.g. \cite{O} if $d=2$ and \cite{B1}, Example 13.1 in general). Then 
$$X'\cong\{([x_0,\ldots,x_d],[y_0,\ldots,y_d])\in\mathbb P^d\times\mathbb P^d\;|\;
x_0y_0+\cdots+x_dy_d=0\}.$$
In case $d=2$ the threefold $X'$ is sometimes called
{\it Hitchin's flag manifold} (see \cite{O}). It follows that
$$\sO_{\mathbb P^d\times\mathbb P^d}(1,1)|X'\cong\sO_{\mathbb P(T_{\mathbb P^d})}(1).$$
Now take $X=\mathbb P^d\times\mathbb P^d$. Then $(\sO_{\mathbb P^d\times\mathbb P^d}(1,1)\cdot Y)_X=
(\sO_{\mathbb P(T_{\mathbb P^d})}(1)\cdot Y)_{X'}=2$. Then Lemma \ref{Lemma3} applies to this situation 
to show that
\begin{equation}\label{eq8}
N_{Y|X}\cong\sO_{\pn 1}(2)\oplus\sO_{\pn 1}(1)^{\oplus 2d-2}.
\end{equation}
In conclusion, $X=\mathbb P^d\times\mathbb P^d$ contains smooth rational curves $Y$ with
the normal bundle given by $(8)$. This example is a higher-dimensional analogue of Example \ref{Example1}.\end{example*}

\medskip

As in the case of surfaces we have the following result:

\begin{prop}\label{Proposition1} Consider the projective variety $X\cong\mathbb P(\sO_{\mathbb P^{n-1}}(2)\oplus\sO_{\mathbb P^{n-1}})$, with $n=2d$ and
$d\geq 2$, and let $Y$ be the smooth rational curve in $X$ with $N_{Y|X}\cong\sO_{\pn 1}(2)\oplus\sO_{\pn 1}(1)^{\oplus 2d-2}$
constructed in Example {\em \ref{Example5}}. Set also $X':=\mathbb P^d\times\mathbb P^d$, and let $Y'$ be the smooth
rational curve in $X'$ with $N_{Y'|X'}\cong\sO_{\pn 1}(2)\oplus\sO_{\pn 1}(1)^{\oplus 2d-2}$
constructed in  Example {\em\ref{Example6}}. Then the pairs $(X,Y)$ and $(X',Y')$ are not formally equivalent.\end{prop}

\proof First we claim that $Y$ is $G3$ in $X$. To see this, clearly we may replace $X$ by the cone $Z$, which is isomorphic to the weighted projective space $\mathbb P^n(1,\ldots,1,2)$. Then the assertion follows  from  \cite{B1}, Corollary 13.3. It also follows that $Y$ meets every hypersurface of $Z$.

On the other hand, since
$N_{Y'|X'}\cong\sO_{\pn 1}(2)\oplus\sO_{\pn 1}(1)^{\oplus 2d-2}$, $N_{Y'|X'}$ is ample, so
by a result of Hartshorne $Y'$ is $G2$ in $X'=\mathbb P^d\times\mathbb P^d$. Since $X'$ is a rational homogeneous space, it follows that $Y'$ is $G3$ in $X$ (see \cite{BSch}, Theorem (4.5), (ii)). Moreover, $Y'$ generates the homogeneous space in the sense of Chow \cite{Ch}; see also \cite{BSch}. Then by Proposition (4.3) of \cite{BSch} it also follows that $Y'$ meets every hypersurface of $X'$. Now assume that the formal completions $X_{/Y}$ and $X'_{/Y'}$ are isomorphic. Then by a  result of Gieseker (see \cite{Gi1} and also \cite{B1}, Corollary 9.20) this implies that there are Zariski open neighbourhoods $U$ in $X$ containing $Y$, and $U'$ in $X'$ containing $Y'$ and a biregular isomorphism $f\colon U\to U'$ such that $f(Y)=Y'$ and $f$ induces  the given formal isomorphism. Again we may replace $X$ by the cone $Z=\mathbb P^n(1,\ldots,1,2)$, which has the advantage that it is a normal $\mathbb Q$-Fano variety. Then the complements $Z\setminus U$ and $X'\setminus U'$ are both of codimension $\geq 2$ (since $Y$ meets every hypersurface of $Z$ and $Y'$ meets every hypersurface of $X'$). The isomorphism $f$ yields an isomorphism between the anticanonical classes $-K_U$ and $-K_{U'}$, and since ${\rm codim}_Z(Z\setminus Y)\geq 2$, $Z$ is a normal $\mathbb Q$-Fano variety, ${\rm codim}_{X'}(X'\setminus Y')\geq 2$ and  $X'$ is a Fano variety, it follows that $f$ extends to an isomorphism 
$Z\cong X'$. But this is absurd because $X'$ is smooth and $Z$ is singular.\qed

\medskip

In dimension $n\geq 3$ there are many more examples of smooth rational curves
$Y$ lying on an $n$-fold $X$ with $N_{Y|X}\cong\sO_{\pn 1}(2)\oplus\sO_{\pn 1}(1)^{\oplus n-2}$ than in dimension $2$, as the following examples show.

\begin{example*}\label{Example7} Let $X'$ be a smooth Fano threefold of index $2$ such that ${\rm Pic}(X')=\mathbb Z[H]$, with $H$ very ample. By a
result of Oxbury \cite{O}  (see also \cite{BBI}, Theorem (3.2), for another proof), $X'$ contains a quasi-line $Y$ which is a conic with respect
to the projective embedding $X'\hookrightarrow\mathbb P^m$ given by $|H|$, i.e. such that $(H\cdot Y)=2$. By Fano-Iskovskih classification
(see \cite{I}), $X'$ is one of the following:
\begin{enumerate}
\item[$-$]  a cubic hypersurface in $\mathbb P^4$ (with $m=4$); or
\item[$-$]  a complete intersection of two hyperquadrics in $\mathbb P^5$ (with $m=5$); or
\item[$-$]  a section of the Pl\"ucker embedding of the Grassmannian $\mathbb G(1,4)$, of lines in $\mathbb P^4$, in $\mathbb P^9$ with three general hyperplanes of $\mathbb P^9$ (with $m=6$).\end{enumerate}

Let now $X$ be a smooth projective fourfold in $\mathbb P^{m+1}$ such that $X'$ is a hyperplane section of $X$. For example, in the first case, 
$X$ can be an arbitrary  cubic fourfold in 
$\mathbb P^5$. Clearly, $N_{X'|X}=\sO_{X'}(1)=H$. Then Lemma \ref{Lemma3} can be applied in this case to get
\begin{equation}\label{eq9}
N_{Y|X}\cong\sO_{\pn 1}(2)\oplus\sO_{\pn 1}(1)\oplus\sO_{\pn 1}(1).\end{equation}
In particular, every cubic fourfold $X$ in $\mathbb P^5$ contains smooth rational curves $Y$ with
the normal bundle given by \eqref{eq9}.\end{example*}

\begin{example*}\label{Example8} Start with the curve $Y=\pn 1\in|\sO(1,1)|$ in $X':=\pn 1\times\mathbb P^1$ of Example \ref{Example1}, and with two
linear embeddings $i\colon\pn 1\hookrightarrow\mathbb P^m$, and
$j\colon\pn 1\hookrightarrow\mathbb P^n$, where $m,n\geq 1$, and $m+n\geq 3$. 
Consider the embedding 
$$i\times j\colon X'=\pn 1\times\mathbb P^1\hookrightarrow X:=\mathbb P^m\times\mathbb P^n.$$
Then $N_{X'|X}\cong\sO_{\mathbb P^m\times\mathbb P^n}(1,0)^{\oplus m-1}\oplus\sO_{\mathbb P^m\times\mathbb P^n}(0,1)^{\oplus n-1}$. In particular,
$N_{X'|X}|Y\cong\sO_{\pn 1}(1)^{\oplus m+n-2}$. Thus the exact sequence
$$0\to N_{Y|X'}=\sO_{\pn 1}(2)\to N_{Y|X}\to N_{X'|X}|Y\cong\sO_{\pn 1}(1)^{\oplus m+n-2}\to 0$$
splits, to give
$$N_{Y|X}\cong\sO_{\pn 1}(2)\oplus\sO_{\pn 1}(1)^{\oplus m+n-2}.$$

If for example we take $m=3$ and $n=1$ we get an embedding $\alpha\colon\pn 1\hookrightarrow\mathbb P^3\times\mathbb P^1$ whose image has 
the normal bundle isomorphic to $\sO_{\pn 1}(2)\oplus\sO_{\pn 1}(1)\oplus\sO_{\pn 1}(1)$. Since the Fano fourfolds
$\mathbb P^2\times\mathbb P^2$ and $\mathbb P^3\times\mathbb P^1$ (which are both homogeneous spaces) cannot be isomorphic
(because $\mathbb P^2\times\mathbb P^2$ has index $3$ and $\mathbb P^3\times\mathbb P^1$ has index $2$), the proof of Proposition \ref{Proposition1} can be applied to yield the fact that this latter embedding cannot be formally equivalent to the embedding $\pn 1\hookrightarrow\mathbb P^2\times\mathbb P^2$ of Example \ref{Example6} (or to the embedding $\beta\colon\pn 1\hookrightarrow\mathbb P^2\times\mathbb P^2$ obtained by the above procedure when $m=n=2$).\end{example*}

\begin{example*}\label{Example9} (Hypercubic in $\mathbb P^{n+1}$) Let $X'$ be a cubic fourfold in $\mathbb P^5$ and $Y\cong\mathbb
P^1\subset X'$ with  $N_{Y|X'}\cong\sO_{\mathbb P^1}(2)\oplus\sO_{\mathbb P^1}(1)\oplus\sO_{\mathbb P^1}(1)$ as in Example \ref{Example7}. Let $X$ be the five
dimensional cubic in $\pn 6$ having $X'$ as hyperplane section. Then $-K_X\cong 4H$ and $X'\in |H|$. Moreover
$N_{X'|X}=H$ and $(Y\cdot H)=2$. Thus arguing as in the proof of Lemma \ref{Lemma3}, we see that the normal bundle of $Y$ in $X$ is
$$N_{Y|X}\cong\sO_{\mathbb P^1}(2)\oplus\sO_{\mathbb P^1}(2)\oplus\sO_{\mathbb P^1}(1)\oplus\sO_{\pn 1}(1).$$
That is, $N_{Y|X}$ is of the form as in Remark  \ref{Remarks}, iii), with
$h:=\dim_\comp L_1(Y,X)=\deg(N_{Y|X})-4=2$.

More generally, we see that an arbitrary hypercubic $X$ in $\mathbb P^{n+1}$ contains a curve $Y\cong\mathbb P^1$ with
$N_{Y|X}\cong\sO_{\mathbb P^1}(2)^{\oplus n-3}\oplus\sO_{\mathbb P^1}(1)\oplus\sO_{\mathbb P^1}(1)$, so that $\deg (N_{Y|X})=2n-4$.
Therefore 
$$h:=\dim_\comp L_1(Y,X)=\deg(N_{Y|X})-n+1=n-3,$$ as in Remark \ref{Remarks}, iii).

Similar conclusions by taking as $(X,H)$ any Del Pezzo $n$-fold with ${\rm Pic}(X)\cong \mathbb Z[H]$. That is $X$ is
the complete intersection of two hyperquadrics in $\mathbb P^{n+2}$, or a linear section of the Grassmannian $\mathbb G(1,4)$ (of lines in $\mathbb P^4$) of dimension $\dim X=4,5$.\end{example*}

\begin{rem*}\label{Fano} Let $X$ be an $n$ dimensional  Fano manifold   containing a smooth rational curve
$Y$  with $N_{Y|X}\cong\sO_{\pn 1}(2)\oplus\sO_{\pn 1}(1)^{\oplus n-2}$. Thus the index, $r$, of $X$ satifies the condition
\begin{equation}\label{index} r\leq\frac{n+2}{2}.\end{equation} 
Indeed, by adjunction formula, $-(K_X\cdot Y)=n+2$. Let $-K_X\cong rH$, $H$ ample line bundle on $X$. Thus $r(H\cdot Y)=n+2$, giving $(H\cdot Y)\geq 2$ and hence the claimed inequality.

In particular, $X$ is neither $\pn n$ nor a hyperquadric. If $X$ is a Del Pezzo manifold (case $r=n-1$) then \eqref{index} yields $n\leq 4$ and therefore $r=3$, $n=4$, as in the case of the cubic fourfold $X$ in $\pn 5$ discussed in Example \ref{Example7}.

\end{rem*}


\bigskip



{\small\begin{thebibliography}{999}

\bibitem{dA}
J. d'Almeida,
Une caract\'erisation du plan projectif complexe, Enseign. Math. \textbf{41} (1995), 135--139.


\bibitem{B1} 
L. B\u adescu,
\textit{Projective Geometry and Formal Geometry},
Monografie Matematyczne Vol. \textbf{65}, Birkh\"auser Verlag (2004).


\bibitem{B2}
L. B\u adescu,
Hyperplane sections and deformations, 
{\em  Algebraic Geometry, Bucharest 1982},
Lecture Notes in Math. \textbf{1054}, Springer-Verlag (1984), pp. 1--33.


\bibitem{BBI}
L. B\u adescu, M.C. Beltrametti and P. Ionescu,
Almost-lines and quasi-lines on projective manifolds, 
\textit{Complex Analysis and Algebraic Geometry - A Volume in Memory
of Michael Schneider}, edited by Peternell and  F.-O. Schreyer, 
Walter de Gruyter, Berlin-New York (2000), pp. 1--27.

\bibitem{BSch}
L. B\u adescu and M. Schneider,
Formal functions, connectivity and homogeneous spaces,
\textit{Algebraic Geometry -- A Volume in Memory of
Paolo Francia}, edited by M.C. Beltrametti, F. Catanese, C. Ciliberto, A. Lanteri and 
C. Pedrini, Walter de Gruyter, Berlin-New York (2002), pp. 1--21.


\bibitem{BS}
M.C. Beltrametti and A.J. Sommese,
\textit{The Adjunction Theory of Complex Projective Varieties}, 
Expositions in Mathematics \textbf{16}, Walter de Gruyter,  Berlin-New York (1995).

\bibitem{Ch}
W.L. Chow,
On meromorphic maps of algebraic varieties,   Annals of Math.
\textbf{89} (1969), 391--403.

\bibitem{De}
M. Demazure,
Surfaces de Del Pezzo,
\textit{S\'eminaire sur les surfaces alg\'ebriques II-V},
Lecture Notes in Math.  \textbf{777} Springer-Verlag (1980), pp. 23--69.

\bibitem{F}
T. Fujita,
On the structure of polarized manifolds with total deficiency one: I,
J. Math. Soc. Japan \textbf{32} (1980), 709--725.


\bibitem{Gi1}
D. Gieseker,
On two theorems of Griffiths about embeddings with ample normal
bundle,  Amer. J. Math.  \textbf{99} (1977),  1137--1150.


\bibitem{Gi2}
D. Gieseker,
$P$-ample bundles and their Chern classes,
Nagoya Math. J.  \textbf{43} (1971), 91--116.


\bibitem{Giz}
M.H. Gizatullin,
On affine surfaces that can be completed by a nonsingular rational curve,
Izv. Akad. Nauk USSR \textbf{34} (1970), 787--810.


\bibitem{GH}
Ph. Griffiths and J. Harris,
\textit{Principles of Algebraic Geometry},
Wiley, Interscience, New York (1978).


\bibitem{FGA}
A. Grothendieck,
\textit{Fondements de la G\'eom\'etrie Alg\'ebrique}, 
Extraits du S\'eminaire Bourbaki,
Paris (1957--1962).


\bibitem{EGA}
A. Grothendieck,
\textit{\'El\'ements de G\'eom\'etrie Alg\'ebrique, III}, Publ. Math. IHES \textbf{11}, Paris (1961).

\bibitem{Ha2}
R. Hartshorne,
\textit{Ample subvarieties of algebraic varieties}, 
Lectures Notes in Math. \textbf{156}, Springer-Verlag, (1970).


\bibitem{Ha1}
 R. Hartshorne,
\textit{Algebraic Geometry}, Graduate Texts in Mathematics, \textbf{52},
 Springer-Verlag (1977).


\bibitem{HM}
H. Hironaka and H. Matsumura,
Formal functions and formal embeddings, 
J. Math. Soc. Japan  \textbf{20} (1968), 52--82.


\bibitem{IN}
P. Ionescu and D. Naie,
Rationality properties of manifolds containing quasi-lines,
Intern. J. Math. \textbf{14} (2003), 1053--1080.

\bibitem{I}
V.A. Iskovskih,
Fano $3$-folds, I,
Math. USSR, Izvestija \textbf{11} (1977), 485--527.


\bibitem{K}
J. Koll\'ar,
\textit{Rational Curves on Algebraic Varieties},
Ergebnisse der Mathematik und ihre Grenzgebiete, \textbf{32},
Springer (1996).


\bibitem{KMM}
J. Koll\'ar, Y. Miyaoka and S. Mori,
Rationally connected varieties,
J. Algebraic Geometry \textbf{1} (1992), 429--448.


\bibitem{MP}
Y. Miyaoka and T. Peternell,
\textit{Geometry of Higher Dimensional Algebraic Varieties},
DMV Seminar \textbf{26}, Birkh\"auser (1997).


\bibitem{O}
W.M. Oxbury,
Twistor spaces and Fano threefolds, 
Quart. J. Math. Oxford \textbf{45} (1994), 343--366.


\bibitem{S}
B. Segre,
\textit{Some Properties of Differentiable Varieties and Transformations -- with Special Reference to 
the Analytic and Algebraic Cases}, Ergebnisse der Mathematik und ihrer Grenzgebiete \textbf{13},
Springer-Verlag (1957).

\end{thebibliography}}
\end{document}